\documentclass[12pt]{amsart}
\usepackage[margin=1.35in]{geometry}

\usepackage{amsmath}
\usepackage{amsthm}
\usepackage{amssymb}
\usepackage{amsfonts}
\usepackage{epsfig}
\usepackage{graphicx,float}
\usepackage{changebar}
\usepackage{amsfonts}
\usepackage{setspace}
\usepackage{mathpazo}
\usepackage{enumerate}
\usepackage{amsthm}
\usepackage{tikz}
\usetikzlibrary{matrix,arrows,decorations.pathmorphing}
\usepackage{lmodern}
\usepackage{calc}
\usepackage{epstopdf}
\usepackage[export]{adjustbox}
\usepackage[hidelinks]{hyperref}
\usepackage{bm}
\usepackage{chngcntr}

\newtheorem{thm}{Theorem}[section]

\newtheorem{lem}[thm]{Lemma}

\newtheorem*{maintheorem}{Main Theorem}
\newenvironment{claim}[1]{\par\noindent\it{Claim:}\space#1}{}

\theoremstyle{definition}
\newtheorem{definition}[thm]{Definition}
\newtheorem{hype}[thm]{Hypothesis}
\newtheorem{example}[thm]{Example}
\newcommand{\bh}{\begin{hyp}\hspace{.225in}}
\newcommand{\eh}{\end{hyp}}
\newcommand{\bt}{\begin{theorem}\hspace{.225in}}
\newcommand{\et}{\end{theorem}}
\newcommand{\br}{\begin{remark}\normalfont\hspace{.225in}}
\newcommand{\er}{\end{remark}}
\newcommand{\bp}{\begin{proposition}\normalfont\hspace{.225in}}
\newcommand{\ep}{\end{proposition}}
\newcommand{\bc}{\begin{corollary}\hspace{.225in}}
\newcommand{\ec}{\end{corollary}}
\newcommand{\bl}{\begin{lemma}\hspace{.225in}}
\newcommand{\el}{\end{lemma}}
\newcommand{\bd}{\begin{definition}\normalfont\hspace{.225in}}
\newcommand{\ed}{\end{definition}}
\newcommand{\be}{\begin{example}\normalfont\hspace{.225in}}
\newcommand{\ee}{\end{example}}
\newcommand{\bpf}{\noindent{\bf Proof.}\indent}
\newcommand{\epf}{\hfill $\blacksquare$\\}
\newcommand{\bcl}{\begin{claim}\hspace{.225in}}
\newcommand{\ecl}{\end{claim}}
\newcommand{\bcpf}{\noindent{\textit {Proof:}}\indent}
\newcommand{\ecpf}{\hfill $\ensuremath{\Box}$}

\numberwithin{equation}{section}

\begin{document}

\title[Classifying character degree graphs with 6 vertices]{Classifying character degree graphs with 6 vertices}

\author[M. W. Bissler]{Mark W. Bissler}
\address{Department of Mathematical Sciences, Kent State University, Kent, OH 44242, U.S.A.}
\email{mbissle2@math.kent.edu}

\author[J. Laubacher]{Jacob Laubacher}
\address{Department of Mathematics, St. Norbert College, De Pere, Wisconsin 54115, U.S.A.}
\email{jacob.laubacher@snc.edu}

\author[M. L. Lewis]{Mark L. Lewis}
\address{Department of Mathematical Sciences, Kent State University, Kent, OH 44242, U.S.A.}
\email{lewis@math.kent.edu}

\subjclass[2010]{20C15, 05C25, 20D10}
\keywords{character degree graphs, solvable groups}
\date{\today}

\begin{abstract}
We investigate prime character degree graphs of solvable groups that have six vertices. There are one hundred twelve non-isomorphic connected graphs with six vertices, of which all except nine are classified in this paper. We also completely classify the disconnected graphs with six vertices.
\end{abstract}

\maketitle

\section{Introduction}

Throughout this paper, $G$ will be a finite solvable group. We will write Irr$(G)$ for the set of irreducible characters of $G$, and cd$(G)=\{\chi(1)\mid \chi\in\text{Irr}(G)\}$. Denote $\rho(G)$ to be the set of primes that divide degrees in cd$(G)$ for the character degrees of $G$. The prime character degree graph of $G$, written $\Delta(G)$, is the graph whose vertex set is $\rho(G)$. Two vertices $p$ and $q$ of $\rho(G)$ are adjacent in $\Delta(G)$ if there exists $a\in\text{cd}(G)$ where $pq$ divides $a$. This type of graph has been studied in a variety of places (see \cite{isaacs}, \cite{lewis2}, \cite{lewis1}, \cite{lewis5}, \cite{palfy}, \cite{zhang}).

In \cite{palfy}, P\'{a}lfy showed that if $G$ is a solvable group for every three vertices in $\rho(G)$, there is some edge in $\Delta(G)$ incident to two of these vertices. With this in mind, we say that a graph $\Gamma$ satisfies P\'{a}lfy's condition if for every three vertices there is some edge incident to two of them. Note that all graphs with three or fewer vertices that satisfy P\'{a}lfy's condition occur as $\Delta(G)$ for some solvable group $G$, as done in \cite{huppert2}. The graphs with four vertices that satisfy P\'{a}lfy's condition all occur except for the disconnected graph with two connected components each of size two and the connected graph with diameter three. These were shown not to occur in \cite{palfy2} and \cite{zhang}, respectively. Lastly, character degree graphs with five vertices were classified in \cite{lewis}, with the exception of one graph.

\begin{figure}[h]\label{Diam3}
    \centering
$
\begin{tikzpicture}[scale=2]
\node (1) at (0,.5) {$\bullet$};
\node (2) at (.5,1) {$\bullet$};
\node (3) at (.5,0) {$\bullet$};
\node (4) at (1,.5) {$\bullet$};
\node (5) at (1.5,.5) {$\bullet$};
\node (6) at (2,.5) {$\bullet$};
\path[font=\small,>=angle 90]
(1) edge node [right] {$ $} (2)
(1) edge node [right] {$ $} (3)
(1) edge node [right] {$ $} (4)
(2) edge node [right] {$ $} (3)
(2) edge node [right] {$ $} (4)
(2) edge node [right] {$ $} (5)
(3) edge node [right] {$ $} (4)
(3) edge node [right] {$ $} (5)
(4) edge node [right] {$ $} (5)
(5) edge node [right] {$ $} (6);
\end{tikzpicture}
$
    \caption{Diameter three graph which occurs as $\Delta(G)$}
    \label{cs}
\end{figure}
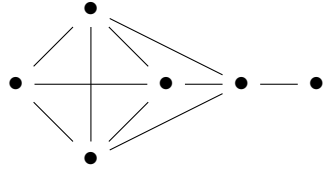

In \cite{lewis4}, the third author settled the conjecture that prime character degree graphs have diameter at most two by constructing the graph in Figure \ref{cs}. Later in \cite{sassy} (see also \cite{sass}), it was shown that this is the unique graph with six vertices and diameter three that occurs as $\Delta(G)$.

Lately, families of graphs have also been studied. In \cite{bissler}, a family of graphs was classified as unable to occur as $\Delta(G)$, and this family was generalized in \cite{laubacher}. These families contain graphs with arbitrarily many vertices, however they do indeed classify several graphs with six vertices along the way, of which we will make use.

In this paper, we look specifically at graphs with six vertices, and we attempt to  determine whether each graph can or cannot occur as the prime character degree graph of a solvable group. In \cite{petric}, it was shown that there are a total of one hundred twelve non-isomorphic connected graphs with six vertices. Of the one hundred twelve graphs, only thirty-five satisfy P\'{a}lfy's condition. We then proceed to ask the following two questions:
\begin{enumerate}
	\item Which of these graphs can occur as the prime character degree graph of some solvable group?
	\item Which of these graphs cannot occur as the prime character degree graph of any solvable group?
\end{enumerate}

After applying known results from the literature and our own results, we are able to conclude that at least thirteen of the connected graphs with six vertices occur as the prime character degree graph for some solvable group, as well as two disconnected graphs. We are left with nine graphs still in consideration. Our results can be summed up into the following theorem:

\begin{maintheorem}
The graphs with six vertices that arise as $\Delta(G)$ for some solvable group $G$ are precisely those graphs in Figures \ref{cs}, \ref{disc}, and \ref{figdirectps}, and possibly those in Figure \ref{figunknown}.
\end{maintheorem}

This work was done while the first author, Mark W. Bissler, was a Ph.D. student at Kent State University under the supervision of Mark L. Lewis. The results in this paper appear in the first author's Ph.D. dissertation \cite{bissler2}. Further work has been conducted with Jacob Laubacher on graphs with six vertices, and has been added.

\section{Preliminaries}

We aim to keep this note as self-contained as possible. We first recall the landmark results from P\'alfy.

\begin{lem}[P\'alfy's condition from \cite{palfy}]
Let $G$ be a solvable group and let $\pi$ be a set of primes contained in $\Delta(G)$. If $|\pi|=3$, then there exists an irreducible character of $G$ with degree divisible by at least two primes from $\pi$. (In other words, any three vertices of the prime character degree graph of a solvable group span at least one edge.)
\end{lem}

In fact, recently P\'alfy's condition was generalized:

\begin{thm}\emph{(\cite{pacifici})}\label{oddc}
Let G be a finite solvable group. Then the complement graph $\bar{\Delta}(G)$ does not contain any cycle of odd length.
\end{thm}

Finally, P\'alfy also has a result which controls the sizes of the components of a disconnected graph. This will be vital when determining which disconnected graphs are possible.

\begin{thm}[P\'alfy's inequality from \cite{palfy2}]\label{palfy2}
Let $G$ be a solvable group and $\Delta(G)$ its prime character degree graph. Suppose that $\Delta(G)$ is disconnected with two components having size $a$ and $b$, where $a\leq b$. Then $b\geq2^a-1$.
\end{thm}

Most of the following results are from \cite{bissler}.

\begin{definition}(\cite{bissler})
A vertex $p$ of a graph $\Gamma$ is \textbf{admissible} if:
\begin{enumerate}[(i)]
\item the subgraph of $\Gamma$ obtained by removing $p$ and all edges incident to $p$ does not occur as the prime character degree graph of any solvable group, and
\item none of the subgraphs of $\Gamma$ obtained by removing one or more of the edges incident to $p$ occur as the prime character degree graph of any solvable group.
\end{enumerate}
\end{definition}

Classifying vertices as admissible is shown to be quite powerful. One of the key results is below, of which we will use often.

\begin{lem}\emph{(\cite{bissler})}
Let G be a solvable group, and suppose $p$ is an admissible vertex of $\Delta(G)$. For every proper normal subgroup $H$ of $G$, suppose that $\Delta(H)$ is a proper subgraph of $\Delta(G)$. Then $O^p(G)=G$.
\end{lem}

En route to showing that a graph does not occur as $\Delta(G)$ for any solvable group $G$, one of the techniques is to first show that $G$ has no normal nonabelian Sylow subgroups.

\begin{lem}\emph{(\cite{bissler})}\label{normiep}
Let $\Gamma$ be a graph satisfying P\'{a}lfy's condition. Let $q$ be a vertex of $\Gamma$, and denote $\pi$ to be the set of vertices of $\Gamma$ adjacent to $q$, and $\rho$ to be the set of vertices not adjacent to $q$. Assume that $\pi$ is the disjoint union of nonempty sets $\pi_1$ and $\pi_2$, and assume that no vertex in $\pi_1$ is adjacent in $\Gamma$ to any vertex in $\pi_2$. Let $v$ be a vertex in $\pi_2$ adjacent to an admissible vertex $s$ in $\rho$. Furthermore, assume there exists another vertex $w$ in $\rho$ that is not adjacent to $v$.

Let $G$ be a solvable group such that $\Delta(G)=\Gamma$, and assume that for every proper normal subgroup $H$ of $G$, $\Delta(H)$ is a proper subgraph of $\Delta(G)$. Then a Sylow $q$-subgroup of $G$ for the prime associated to $q$ is not normal.
\end{lem}
	
We now mention the refinement of admissible which again has consequences involving no normal nonabelian Sylow $p$-subgroup.
    
\begin{definition}(\cite{bissler})
A vertex $p$ of a graph $\Gamma$ is \textbf{strongly admissible} if:
\begin{enumerate}[(i)]
\item $p$ is admissible, and
\item none of the subgraphs of $\Gamma$ obtained by removing $p$, the edges incident to $p$, and one or more of the edges between two adjacent vertices of $p$ occurs as $\Delta(G)$ for some solvable group $G$.
\end{enumerate}
\end{definition}

\begin{lem}\emph{(\cite{bissler})}\label{normp}
Let G be a solvable group, and assume that $p$ is a prime whose vertex is a strongly admissible vertex of $\Delta(G)$. For every proper normal subgroup $H$ of $G$, suppose that $\Delta(G/H)$ is a proper subgraph of $\Delta(G)$. Then a Sylow $p$-subgroup of $G$ is not normal.
\end{lem}

We next mention the result in \cite{lyons}, which is our final method in showing a group $G$ has no normal nonabelian Sylow $p$-subgroup. First, however, we must set some necessary notation. We start by fixing a vertex $p$. Let $\pi$ be the subset of vertices which are adjacent to $p$, and let $\rho$ be the subset of vertices that are not adjacent to $p$. We will let $\pi^*$ and $\rho^*$ denote nonempty subsets of $\pi$ and $\rho$, respectively. Next, for any subset $\pi^*\cup\rho^*$ which induces a complete subgraph of $\Delta(G)$, let $\beta$ denote a subset of vertices which contains $\pi^*\cup\rho^*$ and also induces a complete subgraph of $\Delta(G)$. Setting $\mathcal{B}$ as the union of all such $\beta$'s, we then denote $\tau=\mathcal{B}\setminus(\pi^*\cup\rho^*)$.

\begin{hype}(\cite{lyons})\label{newhype}
Concerning $\Delta(G)$, we assume the following:
\begin{enumerate}[(i)]
    \item\label{00} for every vertex in $\rho$, there exists a nonadjacent vertex in $\pi$,
    \item\label{11} for every vertex in $\pi$, there exists a nonadjacent vertex in $\rho$,
    \item\label{33} all the vertices in $\pi$ are admissible. Moreover, no proper connected subgraph with vertex set $\{p\}\cup\pi^*\cup\rho$ occurs as the prime character degree graph of any solvable group,
    \item\label{44} for each vertex set $\pi^*\cup\rho^*$ which induces a complete subgraph in $\Delta(G)$, all the vertices in the corresponding set $\tau$ are admissible. Moreover, no proper connected subgraph with vertex set $\rho(G)\setminus\tau^*$ occurs as the prime character degree graph of any solvable group, and
    \item\label{55} if a disconnected subgraph with vertex set $\rho(G)$ does not occur, then it must specifically violate P\'alfy's inequality from \cite{palfy2}. Finally, if a disconnected subgraph with vertex set $\rho(G)$ does occur, then the sizes of the connected components must be $n>1$ and $2^n-1$.
\end{enumerate}
\end{hype}

\begin{thm}\emph{(\cite{lyons})}\label{hypethm}
Assume Hypothesis \ref{newhype}. Then $G$ has no normal nonabelian Sylow $p$-subgroup.
\end{thm}  
 
The penultimate result we provide in this section will be the tool we use to show a graph cannot occur for any solvable group $G$. Note this is Lemma 2.6 of \cite{bissler}.

\begin{lem}\emph{(\cite{bissler})}\label{final}
Let $\Gamma$ be a graph satisfying P\'alfy's condition with $n\geq5$ vertices. Also, assume there exist distinct vertices $a$ and $b$ of $\Gamma$ such that $a$ is adjacent to an admissible vertex $c$, $b$ is not adjacent to $c$, and $a$ is not adjacent to an admissible vertex $d$.

Let $G$ be a solvable group and suppose for all proper normal subgroups $N$ of $G$ we have that $\Delta(N)$ and $\Delta(G/N)$ are proper subgraphs of $\Gamma$. Let $F$ be the Fitting subgroup of $G$ and suppose that $F$ in minimal normal in $G$. Then $\Gamma$ is not the prime character degree graph of any solvable group.
\end{lem}

Finally, we recall the main result from \cite{bissler}. This constructs a family of graphs which cannot occur as $\Delta(G)$ for any solvable group $G$. This will be useful since the family contains graphs with six vertices.

\begin{thm}\emph{(\cite{bissler})}\label{bissmain}
Let $\Gamma$ be a graph satisfying P\'alfy's condition with $k\geq5$ vertices. Assume that there exists two vertices $p_1$ and $p_2$ in $\Gamma$ such that
\begin{enumerate}[(i)]
    \item both $p_1$ and $p_2$ are of degree two,
    \item $p_1$ is adjacent to $p_2$, and
    \item $p_1$ and $p_2$ share no common neighbor.
\end{enumerate}
Then $\Gamma$ is not the prime character degree graph of any solvable group.
\end{thm}

\section{Constructions}\label{const}

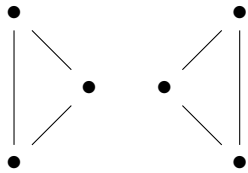
\begin{figure}[h]
    \centering
$
\begin{tikzpicture}[scale=2]
\node (a) at (.5,1) {$\bullet$};
\node (b) at (.5,0) {$\bullet$};
\node (c) at (1,.5) {$\bullet$};
\node (d) at (2,1) {$\bullet$};
\node (e) at (2,0) {$\bullet$};
\node (f) at (1.5,.5) {$\bullet$};
\path[font=\small,>=angle 90]
(a) edge node [right] {$ $} (b)
(a) edge node [above] {$ $} (c)
(b) edge node [above] {$ $} (c)
(d) edge node [above] {$ $} (e)
(d) edge node [above] {$ $} (f)
(e) edge node [right] {$ $} (f);
\end{tikzpicture}
$
    \caption{Disconnected graph which does not occur as $\Delta(G)$}
    \label{fig3and3}
\end{figure}

Be begin this section by considering the disconnected graphs with six vertices. Due to P\'alfy's condition, we can have at most two components. Moreover, each component must be a complete graph. Considering graphs with six vertices, therefore, yields three possibilities: graphs with component sizes one and five, two and four, and three and three. We know the graph with component sizes three and three cannot occur (see Figure \ref{fig3and3}). This is due to P\'alfy's inequality by taking $a=b=3$ and observing that $3\geq7$, a contradiction.

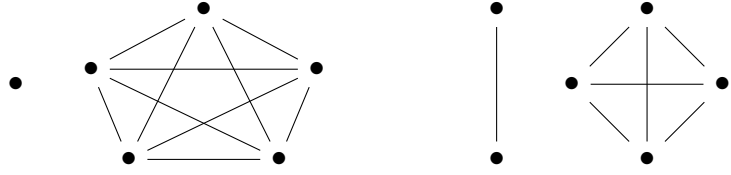
\begin{figure}[t]
    \centering
$
\begin{tikzpicture}[scale=2]
\node (i) at (0,.6) {$\bullet$};
\node (a) at (.25,0) {$\bullet$};
\node (b) at (1.25,0) {$\bullet$};
\node (c) at (1.5,.6) {$\bullet$};
\node (d) at (.75,1) {$\bullet$};
\node (f) at (-.5,.5) {$\bullet$};
\path[font=\small,>=angle 90]
(i) edge node [right] {$ $} (a)
(i) edge node [above] {$ $} (b)
(i) edge node [above] {$ $} (c)
(i) edge node [above] {$ $} (d)
(a) edge node [right] {$ $} (b)
(a) edge node [above] {$ $} (c)
(a) edge node [above] {$ $} (d)
(b) edge node [above] {$ $} (c)
(b) edge node [above] {$ $} (d)
(c) edge node [right] {$ $} (d);
\end{tikzpicture}
\hspace{.75in}
\begin{tikzpicture}[scale=2]
\node (a) at (0,.5) {$\bullet$};
\node (b) at (.5,0) {$\bullet$};
\node (c) at (.5,1) {$\bullet$};
\node (d) at (1,.5) {$\bullet$};
\node (f) at (-.5,0) {$\bullet$};
\node (g) at (-.5,1) {$\bullet$};
\path[font=\small,>=angle 90]
(a) edge node [right] {$ $} (b)
(a) edge node [above] {$ $} (c)
(a) edge node [above] {$ $} (d)
(b) edge node [above] {$ $} (c)
(b) edge node [above] {$ $} (d)
(c) edge node [right] {$ $} (d)
(f) edge node [above] {$ $} (g);
\end{tikzpicture}
$
    \caption{Disconnected graphs which occur as $\Delta(G)$}
    \label{disc}
\end{figure}

We provide the construction of the two possible disconnected graphs, displayed in Figure \ref{disc}. We mimic the process for the construction of these from \cite{lewis}. First, we consider the construction of the graph with component sizes one and five. For this, consider the field of order $2^{32}$ acted on by its full multiplication group and then its Galois group. This group has character degree set $\{1,2,4,8,16,32,3\cdot 5\cdot 17\cdot 257\cdot 65537 \}$. Thus, this gives a graph with two connected components: $\{2\}$ and $\{3,5,17,257,65537\}$.

Next, we consider the graph with component sizes two and four. Consider the field of order $2^{35}$ acted on by its full multiplication group and then its Galois group. This group has character degree set $\{1,5,7,35,31\cdot 71\cdot 127\cdot 122921\}$. Thus, this gives a graph with two connected components: $\{5,7\}$ and $\{31,71,127,122921\}$.

We also construct the graphs that can occur as prime character degree graphs via direct products. For groups $H$ and $K$ we have that $\rho(H\times K)=\rho(H)\cup\rho(K)$. We require also that $\rho(H)$ and $\rho(K)$ be disjoint sets; that is they must have different primes. Since there are infinitely many prime sets that occur for groups having these graphs, it is not a problem to find such groups $H$ and $K$. Next, there is an edge between vertices $p$ and $q$ in $\Delta(H\times K)$ if we have any of the following cases:
\begin{enumerate}
\item $p,q\in\rho(H)$ and there is an edge between $p$ and $q$ in $\Delta(H)$, 
\item $p,q\in\rho(K)$ and there is an edge between $p$ and $q$ in $\Delta(K)$, 
\item $p\in\rho(H)$ and $q\in\rho(K)$, or 
\item $p\in\rho(K)$ and $q\in\rho(H)$. 
\end{enumerate}

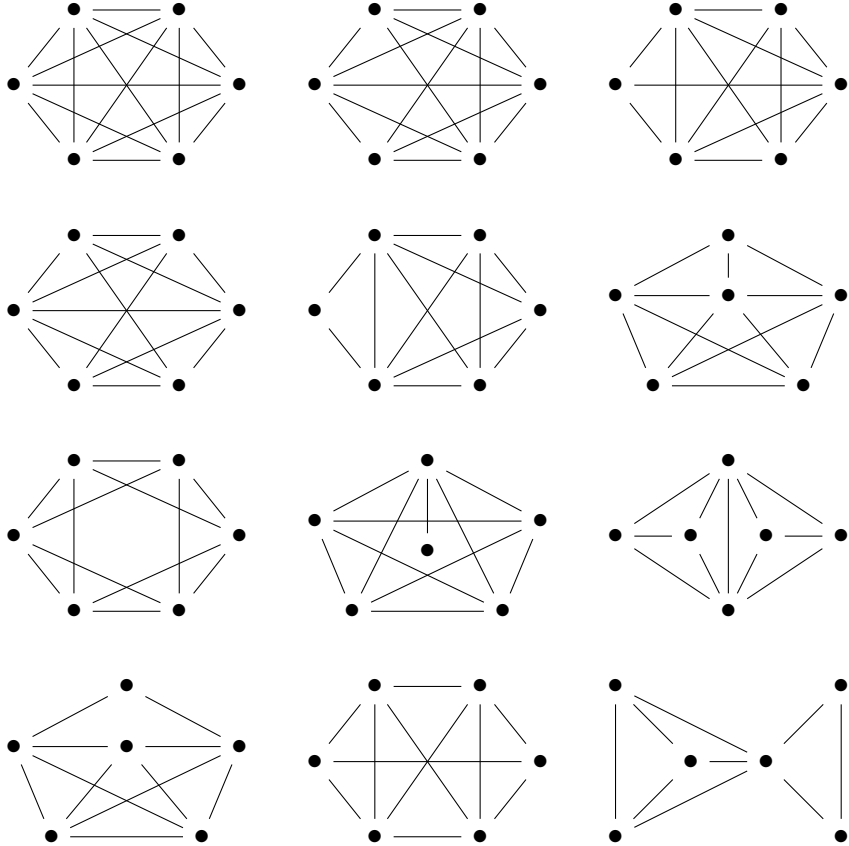
\begin{figure}[h]\label{directps}
    \centering
$
\begin{tikzpicture}[scale=2]
\node (1a) at (0,5) {$\bullet$};
\node (1b) at (.4,5.5) {$\bullet$};
\node (1c) at (.4,4.5) {$\bullet$};
\node (1d) at (1.1,5.5) {$\bullet$};
\node (1e) at (1.1,4.5) {$\bullet$};
\node (1f) at (1.5,5) {$\bullet$};
\path[font=\small,>=angle 90]
(1a) edge node [right] {$ $} (1b)
(1a) edge node [right] {$ $} (1c)
(1a) edge node [right] {$ $} (1d)
(1a) edge node [right] {$ $} (1e)
(1a) edge node [right] {$ $} (1f)
(1b) edge node [right] {$ $} (1c)
(1b) edge node [right] {$ $} (1d)
(1b) edge node [right] {$ $} (1e)
(1b) edge node [right] {$ $} (1f)
(1c) edge node [right] {$ $} (1d)
(1c) edge node [right] {$ $} (1e)
(1c) edge node [right] {$ $} (1f)
(1d) edge node [right] {$ $} (1e)
(1d) edge node [right] {$ $} (1f)
(1e) edge node [right] {$ $} (1f);
\node (2a) at (2,5) {$\bullet$};
\node (2b) at (2.4,5.5) {$\bullet$};
\node (2c) at (2.4,4.5) {$\bullet$};
\node (2d) at (3.1,5.5) {$\bullet$};
\node (2e) at (3.1,4.5) {$\bullet$};
\node (2f) at (3.5,5) {$\bullet$};
\path[font=\small,>=angle 90]
(2a) edge node [right] {$ $} (2b)
(2a) edge node [right] {$ $} (2c)
(2a) edge node [right] {$ $} (2d)
(2a) edge node [right] {$ $} (2e)
(2a) edge node [right] {$ $} (2f)
(2b) edge node [right] {$ $} (2d)
(2b) edge node [right] {$ $} (2e)
(2b) edge node [right] {$ $} (2f)
(2c) edge node [right] {$ $} (2d)
(2c) edge node [right] {$ $} (2e)
(2c) edge node [right] {$ $} (2f)
(2d) edge node [right] {$ $} (2e)
(2d) edge node [right] {$ $} (2f)
(2e) edge node [right] {$ $} (2f);
\node (3a) at (4,5) {$\bullet$};
\node (3b) at (4.4,5.5) {$\bullet$};
\node (3c) at (4.4,4.5) {$\bullet$};
\node (3d) at (5.1,5.5) {$\bullet$};
\node (3e) at (5.1,4.5) {$\bullet$};
\node (3f) at (5.5,5) {$\bullet$};
\path[font=\small,>=angle 90]
(3a) edge node [right] {$ $} (3b)
(3a) edge node [right] {$ $} (3c)
(3a) edge node [right] {$ $} (3f)
(3b) edge node [right] {$ $} (3c)
(3b) edge node [right] {$ $} (3d)
(3b) edge node [right] {$ $} (3e)
(3b) edge node [right] {$ $} (3f)
(3c) edge node [right] {$ $} (3d)
(3c) edge node [right] {$ $} (3e)
(3c) edge node [right] {$ $} (3f)
(3d) edge node [right] {$ $} (3e)
(3d) edge node [right] {$ $} (3f)
(3e) edge node [right] {$ $} (3f);
\node (4a) at (0,3.5) {$\bullet$};
\node (4b) at (.4,4) {$\bullet$};
\node (4c) at (.4,3) {$\bullet$};
\node (4d) at (1.1,4) {$\bullet$};
\node (4e) at (1.1,3) {$\bullet$};
\node (4f) at (1.5,3.5) {$\bullet$};
\path[font=\small,>=angle 90]
(4a) edge node [right] {$ $} (4b)
(4a) edge node [right] {$ $} (4c)
(4a) edge node [right] {$ $} (4d)
(4a) edge node [right] {$ $} (4e)
(4a) edge node [right] {$ $} (4f)
(4b) edge node [right] {$ $} (4d)
(4b) edge node [right] {$ $} (4e)
(4b) edge node [right] {$ $} (4f)
(4c) edge node [right] {$ $} (4d)
(4c) edge node [right] {$ $} (4e)
(4c) edge node [right] {$ $} (4f)
(4d) edge node [right] {$ $} (4f)
(4e) edge node [right] {$ $} (4f);
\node (5a) at (2,3.5) {$\bullet$};
\node (5b) at (2.4,4) {$\bullet$};
\node (5c) at (2.4,3) {$\bullet$};
\node (5d) at (3.1,4) {$\bullet$};
\node (5e) at (3.1,3) {$\bullet$};
\node (5f) at (3.5,3.5) {$\bullet$};
\path[font=\small,>=angle 90]
(5a) edge node [right] {$ $} (5b)
(5a) edge node [right] {$ $} (5c)
(5b) edge node [right] {$ $} (5c)
(5b) edge node [right] {$ $} (5d)
(5b) edge node [right] {$ $} (5e)
(5b) edge node [right] {$ $} (5f)
(5c) edge node [right] {$ $} (5d)
(5c) edge node [right] {$ $} (5e)
(5c) edge node [right] {$ $} (5f)
(5d) edge node [right] {$ $} (5e)
(5d) edge node [right] {$ $} (5f)
(5e) edge node [right] {$ $} (5f);
\node (6a) at (4.75,4) {$\bullet$};
\node (6b) at (4,3.6) {$\bullet$};
\node (6c) at (4.75,3.6) {$\bullet$};
\node (6d) at (5.5,3.6) {$\bullet$};
\node (6e) at (4.25,3) {$\bullet$};
\node (6f) at (5.25,3) {$\bullet$};
\path[font=\small,>=angle 90]
(6a) edge node [right] {$ $} (6b)
(6a) edge node [right] {$ $} (6c)
(6a) edge node [right] {$ $} (6d)
(6b) edge node [right] {$ $} (6c)
(6b) edge node [right] {$ $} (6e)
(6b) edge node [right] {$ $} (6f)
(6c) edge node [right] {$ $} (6d)
(6c) edge node [right] {$ $} (6e)
(6c) edge node [right] {$ $} (6f)
(6d) edge node [right] {$ $} (6e)
(6d) edge node [right] {$ $} (6f)
(6e) edge node [right] {$ $} (6f);
\node (7a) at (0,2) {$\bullet$};
\node (7b) at (.4,2.5) {$\bullet$};
\node (7c) at (.4,1.5) {$\bullet$};
\node (7d) at (1.1,2.5) {$\bullet$};
\node (7e) at (1.1,1.5) {$\bullet$};
\node (7f) at (1.5,2) {$\bullet$};
\path[font=\small,>=angle 90]
(7a) edge node [right] {$ $} (7b)
(7a) edge node [right] {$ $} (7c)
(7a) edge node [right] {$ $} (7d)
(7a) edge node [right] {$ $} (7e)
(7b) edge node [right] {$ $} (7c)
(7b) edge node [right] {$ $} (7d)
(7b) edge node [right] {$ $} (7f)
(7c) edge node [right] {$ $} (7e)
(7c) edge node [right] {$ $} (7f)
(7d) edge node [right] {$ $} (7e)
(7d) edge node [right] {$ $} (7f)
(7e) edge node [right] {$ $} (7f);
\node (8a) at (2.75,1.9) {$\bullet$};
\node (8b) at (2,2.1) {$\bullet$};
\node (8c) at (2.75,2.5) {$\bullet$};
\node (8d) at (3.5,2.1) {$\bullet$};
\node (8e) at (2.25,1.5) {$\bullet$};
\node (8f) at (3.25,1.5) {$\bullet$};
\path[font=\small,>=angle 90]
(8a) edge node [right] {$ $} (8c)
(8b) edge node [right] {$ $} (8c)
(8b) edge node [right] {$ $} (8d)
(8b) edge node [right] {$ $} (8e)
(8b) edge node [right] {$ $} (8f)
(8c) edge node [right] {$ $} (8d)
(8c) edge node [right] {$ $} (8e)
(8c) edge node [right] {$ $} (8f)
(8d) edge node [right] {$ $} (8e)
(8d) edge node [right] {$ $} (8f)
(8e) edge node [right] {$ $} (8f);
\node (9a) at (4,2) {$\bullet$};
\node (9b) at (4.5,2) {$\bullet$};
\node (9c) at (4.75,2.5) {$\bullet$};
\node (9d) at (4.75,1.5) {$\bullet$};
\node (9e) at (5,2) {$\bullet$};
\node (9f) at (5.5,2) {$\bullet$};
\path[font=\small,>=angle 90]
(9a) edge node [right] {$ $} (9b)
(9a) edge node [right] {$ $} (9c)
(9a) edge node [right] {$ $} (9d)
(9b) edge node [right] {$ $} (9c)
(9b) edge node [right] {$ $} (9d)
(9c) edge node [right] {$ $} (9d)
(9c) edge node [right] {$ $} (9e)
(9c) edge node [right] {$ $} (9f)
(9d) edge node [right] {$ $} (9e)
(9d) edge node [right] {$ $} (9f)
(9e) edge node [right] {$ $} (9f);
\node (10a) at (.75,1) {$\bullet$};
\node (10b) at (0,.6) {$\bullet$};
\node (10c) at (.75,.6) {$\bullet$};
\node (10d) at (1.5,.6) {$\bullet$};
\node (10e) at (.25,0) {$\bullet$};
\node (10f) at (1.25,0) {$\bullet$};
\path[font=\small,>=angle 90]
(10a) edge node [right] {$ $} (10b)
(10a) edge node [right] {$ $} (10d)
(10b) edge node [right] {$ $} (10c)
(10b) edge node [right] {$ $} (10e)
(10b) edge node [right] {$ $} (10f)
(10c) edge node [right] {$ $} (10d)
(10c) edge node [right] {$ $} (10e)
(10c) edge node [right] {$ $} (10f)
(10d) edge node [right] {$ $} (10e)
(10d) edge node [right] {$ $} (10f)
(10e) edge node [right] {$ $} (10f);
\node (11a) at (2,.5) {$\bullet$};
\node (11b) at (2.4,1) {$\bullet$};
\node (11c) at (2.4,0) {$\bullet$};
\node (11d) at (3.1,1) {$\bullet$};
\node (11e) at (3.1,0) {$\bullet$};
\node (11f) at (3.5,.5) {$\bullet$};
\path[font=\small,>=angle 90]
(11a) edge node [right] {$ $} (11b)
(11a) edge node [right] {$ $} (11c)
(11a) edge node [right] {$ $} (11f)
(11b) edge node [right] {$ $} (11c)
(11b) edge node [right] {$ $} (11d)
(11b) edge node [right] {$ $} (11e)
(11c) edge node [right] {$ $} (11d)
(11c) edge node [right] {$ $} (11e)
(11d) edge node [right] {$ $} (11e)
(11d) edge node [right] {$ $} (11f)
(11e) edge node [right] {$ $} (11f);
\node (12a) at (4,1) {$\bullet$};
\node (12b) at (4,0) {$\bullet$};
\node (12c) at (4.5,.5) {$\bullet$};
\node (12d) at (5,.5) {$\bullet$};
\node (12e) at (5.5,1) {$\bullet$};
\node (12f) at (5.5,0) {$\bullet$};
\path[font=\small,>=angle 90]
(12a) edge node [right] {$ $} (12b)
(12a) edge node [right] {$ $} (12c)
(12a) edge node [right] {$ $} (12d)
(12b) edge node [right] {$ $} (12c)
(12b) edge node [right] {$ $} (12d)
(12c) edge node [right] {$ $} (12d)
(12d) edge node [right] {$ $} (12e)
(12d) edge node [right] {$ $} (12f)
(12e) edge node [right] {$ $} (12f);
\end{tikzpicture}
$
    \caption{Direct products which occur as $\Delta(G)$}
    \label{figdirectps}
\end{figure}

	Using this construction, we are able to produce the twelve graphs in Figure \ref{figdirectps} that occur as $\Delta(G)$ for some solvable group $G$.

\section{Reductions}

We start by noting the number of connected graphs with six vertices is one hundred twelve by \cite{petric}. After applying P\'{a}lfy's condition to the one hundred twelve graphs, we reduce this number to thirty-five possible graphs. We use the results of \cite{sass} when dealing with a graph or subgraph that arises with diameter three. These results from \cite{sass} (see also her dissertation \cite{sassy}) conclude that the graph in Figure \ref{cs} is the unique six-vertex graph with diameter three that occurs as $\Delta(G)$. Using this fact we can decrease the number of graphs from thirty-five graphs down to twenty-seven.

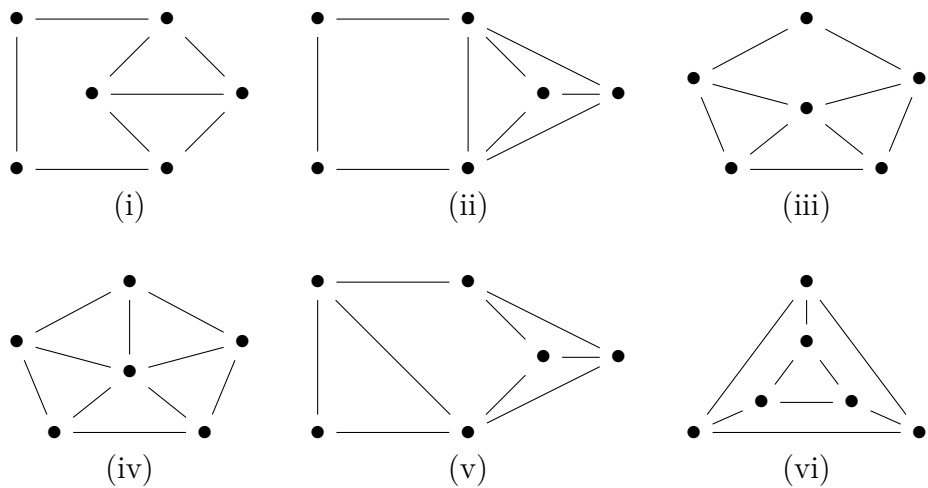
\begin{figure}[h]
    \centering
$
\begin{tikzpicture}[scale=2]
\node (3a) at (0,1) {$\bullet$};
\node (3b) at (0,0) {$\bullet$};
\node (3c) at (1,1) {$\bullet$};
\node (3d) at (1,0) {$\bullet$};
\node (3e) at (.5,.5) {$\bullet$};
\node (3f) at (1.5,.5) {$\bullet$};
\node (3z) at (.75,-.25) {(i)};
\path[font=\small,>=angle 90]
(3a) edge node [right] {$ $} (3b)
(3a) edge node [right] {$ $} (3c)
(3b) edge node [right] {$ $} (3d)
(3c) edge node [right] {$ $} (3e)
(3c) edge node [right] {$ $} (3f)
(3d) edge node [right] {$ $} (3e)
(3d) edge node [right] {$ $} (3f)
(3e) edge node [right] {$ $} (3f);
\node (1a) at (2,1) {$\bullet$};
\node (1b) at (2,0) {$\bullet$};
\node (1c) at (3,1) {$\bullet$};
\node (1d) at (3,0) {$\bullet$};
\node (1e) at (3.5,.5) {$\bullet$};
\node (1f) at (4,.5) {$\bullet$};
\node (1z) at (3,-.25) {(ii)};
\path[font=\small,>=angle 90]
(1a) edge node [right] {$ $} (1b)
(1a) edge node [right] {$ $} (1c)
(1b) edge node [right] {$ $} (1d)
(1c) edge node [right] {$ $} (1d)
(1c) edge node [right] {$ $} (1e)
(1c) edge node [right] {$ $} (1f)
(1d) edge node [right] {$ $} (1e)
(1d) edge node [right] {$ $} (1f)
(1e) edge node [right] {$ $} (1f);
\node (4a) at (5.25,1) {$\bullet$};
\node (4b) at (4.5,.6) {$\bullet$};
\node (4c) at (5.25,.4) {$\bullet$};
\node (4d) at (6,.6) {$\bullet$};
\node (4e) at (4.75,0) {$\bullet$};
\node (4f) at (5.75,0) {$\bullet$};
\node (4z) at (5.25,-.25) {(iii)};
\path[font=\small,>=angle 90]
(4a) edge node [right] {$ $} (4b)
(4a) edge node [right] {$ $} (4d)
(4b) edge node [right] {$ $} (4c)
(4b) edge node [right] {$ $} (4e)
(4c) edge node [right] {$ $} (4d)
(4c) edge node [right] {$ $} (4e)
(4c) edge node [right] {$ $} (4f)
(4d) edge node [right] {$ $} (4f)
(4e) edge node [right] {$ $} (4f);
\node (i) at (0,-1.15) {$\bullet$};
\node (a) at (.25,-1.75) {$\bullet$};
\node (b) at (1.25,-1.75) {$\bullet$};
\node (c) at (1.5,-1.15) {$\bullet$};
\node (d) at (.75,-.75) {$\bullet$};
\node (f) at (.75,-1.35) {$\bullet$};
\node (i0) at (.75,-2) {(iv)};
\path[font=\small,>=angle 90]
(f) edge node [right] {$ $} (i)
(f) edge node [above] {$ $} (a)
(f) edge node [above] {$ $} (b)
(f) edge node [above] {$ $} (c)
(f) edge node [right] {$ $} (d)
(i) edge node [above] {$ $} (a)
(a) edge node [above] {$ $} (b)
(b) edge node [above] {$ $} (c)
(c) edge node [above] {$ $} (d)
(d) edge node [right] {$ $} (i);
\node (2a) at (2,-.75) {$\bullet$};
\node (2b) at (2,-1.75) {$\bullet$};
\node (2c) at (3,-.75) {$\bullet$};
\node (2d) at (3,-1.75) {$\bullet$};
\node (2e) at (3.5,-1.25) {$\bullet$};
\node (2f) at (4,-1.25) {$\bullet$};
\node (2z) at (3,-2) {(v)};
\path[font=\small,>=angle 90]
(2a) edge node [right] {$ $} (2b)
(2a) edge node [right] {$ $} (2c)
(2a) edge node [right] {$ $} (2d)
(2b) edge node [right] {$ $} (2d)
(2c) edge node [right] {$ $} (2e)
(2c) edge node [right] {$ $} (2f)
(2d) edge node [right] {$ $} (2e)
(2d) edge node [right] {$ $} (2f)
(2e) edge node [right] {$ $} (2f);
\node (5a) at (5.25,-.75) {$\bullet$};
\node (5b) at (4.5,-1.75) {$\bullet$};
\node (5c) at (6,-1.75) {$\bullet$};
\node (5d) at (5.25,-1.15) {$\bullet$};
\node (5e) at (4.95,-1.55) {$\bullet$};
\node (5f) at (5.55,-1.55) {$\bullet$};
\node (5z) at (5.25,-2) {(vi)};
\path[font=\small,>=angle 90]
(5a) edge node [right] {$ $} (5b)
(5a) edge node [right] {$ $} (5c)
(5a) edge node [right] {$ $} (5d)
(5b) edge node [right] {$ $} (5c)
(5b) edge node [right] {$ $} (5e)
(5c) edge node [right] {$ $} (5f)
(5d) edge node [right] {$ $} (5e)
(5d) edge node [right] {$ $} (5f)
(5e) edge node [right] {$ $} (5f);
\end{tikzpicture}
$
    \caption{Connected six-vertex graphs investigated in this paper}
    \label{fig6graphs}
\end{figure}

	Recall that in Section \ref{const} we constructed twelve graphs via direct products (see Figure \ref{directps}); if we remove these twelve graphs from the reduced list of twenty-seven, we are then left with fifteen graphs. In this section we will show that the graphs in Figure \ref{fig6graphs} cannot occur as $\Delta(G)$ for any solvable group $G$.

\begin{thm}
Graphs (i) and (ii) of Figure \ref{fig6graphs} are not the prime character degree graph of a solvable group.
\end{thm}
  
\bpf
Using the main result from \cite{bissler} (see Theorem \ref{bissmain} above), we have that graphs (i) and (ii) from Figure \ref{fig6graphs} fall into the infinite family shown not to occur as the prime character degree graph for any solvable group. One observes that these are also classified in the families studied in \cite{laubacher}.
\epf 

\begin{thm}
Graphs (iii) and (iv) of Figure \ref{fig6graphs} are not the prime character degree graph of a solvable group.
\end{thm}

\bpf 
By the main result from \cite{pacifici} (see Theorem \ref{oddc} above), any graph whose complement graph has an odd cycle can be eliminated. The complement of graphs (iii) and (iv) contain a cycle of length five, thus they cannot be the prime character degree graph for a solvable group.
\epf

\begin{thm}
Graph (vi) of Figure \ref{fig6graphs} is not the prime character degree graph of a solvable group.
\end{thm}

\bpf 
In \cite{laubacher}, the authors show that the graph (vi) of Figure \ref{fig6graphs} cannot occur as $\Delta(G)$ for any solvable group $G$ by showing that every vertex is admissible. This graph is also classified in an infinite family which cannot occur as the character degree graph for any solvable group. We mention also, that in \cite{Paolo}, a result about regular character degree graphs was shown which also proves graph (vi) is not the prime character degree graph of a solvable group.
\epf

\subsection{Eliminating graph (v) of Figure \ref{fig6graphs}}

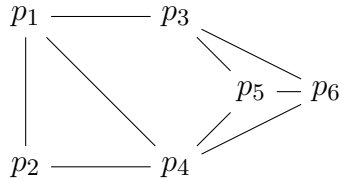
\begin{figure}[h]
\centering
$
\begin{tikzpicture}[scale=2]
\node (1a) at (2,1) {$p_1$};
\node (1b) at (2,0) {$p_2$};
\node (1c) at (3,1) {$p_3$};
\node (1d) at (3,0) {$p_4$};
\node (1e) at (3.5,.5) {$p_5$};
\node (1f) at (4,.5) {$p_6$};
\path[font=\small,>=angle 90]
(1a) edge node [right] {$ $} (1b)
(1a) edge node [right] {$ $} (1c)
(1b) edge node [right] {$ $} (1d)
(1a) edge node [right] {$ $} (1d)
(1c) edge node [right] {$ $} (1e)
(1c) edge node [right] {$ $} (1f)
(1d) edge node [right] {$ $} (1e)
(1d) edge node [right] {$ $} (1f)
(1e) edge node [right] {$ $} (1f);
\end{tikzpicture}
$
\caption{The graph in Theorem \ref{th3}}
\label{figf}
\end{figure}

Due to the length of the proof, we dedicate this separate subsection to eliminating graph (v) of Figure \ref{fig6graphs}. It is beneficial to label our vertices, as in Figure \ref{figf}.

\begin{thm}\label{th3}
There is no solvable group $G$ such that $\Delta(G)$ is the graph found in Figure \ref{figf}.
\end{thm}

\bpf
Suppose $G$ is a counterexample with $|G|$ minimal. We label the primes in $\rho(G)$ as the vertices in Figure \ref{figf}. Our first goal will be to show that $G$ does not have any normal nonabelian Sylow subgroups. 
We first show which vertices of Figure \ref{figf} are strongly admissible.\\

\bcl
The vertices $p_1$, $p_4$, $p_5$, and $p_6$ are strongly admissible vertices, and thus $G$ does not have a corresponding normal Sylow subgroup for $p_1$, $p_4$, $p_5$ and $p_6$.
\ecl\\

\bcpf
We start by considering the vertex $p_1$. If we lose an edge between $p_1$ and $p_3$, then we arrive at a graph with diameter three, which is not the unique one, so we cannot lose this edge. If we lose an edge between $p_1$ and $p_2$, we violate P\'{a}lfy's condition with $p_1$, $p_2$, and $p_6$. If we lose the edge between $p_1$ and $p_4$, we have the graph (i) of Figure \ref{fig6graphs}, which has been shown to not occur. Now let us assume that the vertex $p_1$ is removed from the graph. Let $P_1$ be a Sylow $p_1$-subgroup of $G$ and assume that $P_1\triangleleft G$. By Lemma 3 of \cite{lewis1} we know $\mbox{cd}(G/P_1')=\mbox{cd}(G)\setminus\{p_1\}$, so the subgraph will have vertex set $\{p_2, p_3, p_4, p_5, p_6\}$. If no edges are lost, we arrive at a subgraph with five vertices that has diameter three, which we know is not possible by the main theorem of \cite{lewis2}. The only possible edge that could be lost in this subgraph is the edge between $p_2$ and $p_4$, and this subgraph would violate P\'{a}lfy's condition. Thus, $p_1$ is strongly admissible.

Now we look at the vertex $p_4$. If an edge is lost between $p_2$ and $p_4$, we violate P\'{a}lfy's condition with vertices $p_2$, $p_3$, and $p_4$. We have already shown we cannot lose the edge between $p_1$ and $p_4$, so we consider the edge between $p_4$ and $p_5$ or $p_6$. If the edge between $p_4$ and $p_5$ is removed, we then arrive at a graph with diameter three, which is not the unique diameter three graph. If we lose the edge between $p_4$ and $p_6$, we then also arrive at a diameter three graph that is not the unique diameter three graph with six vertices. Let $P_4$ be a Sylow $p_4$-subgroup of $G$ and assume that $P_4\triangleleft G$. By Lemma 3 of \cite{lewis1} we know $\mbox{cd}(G/P_4')=\mbox{cd}(G)\setminus\{p_1\}$, so the subgraph will have vertex set $\{p_1, p_2, p_3, p_5, p_6\}$. If we remove the vertex $p_4$ from the graph, the subgraph that loses no other edges has fives vertices and diameter three, which we know is not possible again, by the main theorem of \cite{lewis2}. The only edges in this subgraph that could be lost are the edges between $p_1$ and $p_2$ and the edge between $p_5$ and $p_6$. In either case we result in a subgraph that violates P\'{a}lfy's condition. Thus, $p_4$ is strongly admissible. 

We next show $p_5$ is strongly admissible and note that a similar argument will show that $p_6$ is strongly admissible. If an edge is lost between $p_3$ and $p_5$, then we violate P\'{a}lfy's condition with $p_3$, $p_5$, and $p_2$. If an edge is lost between $p_5$ and $p_6$, we again violate P\'{a}lfy's condition with $p_5$, $p_6$, and $p_2$. The edge between $p_4$ and $p_5$ we have already shown cannot be lost, so we now consider when the vertex  $p_5$ is removed from the graph. Let $P_5$ be a Sylow $p_5$-subgroup of $G$ and assume that $P_5\triangleleft G$. If $p_5$ is removed from the graph, then by Lemma 3 of \cite{lewis1}, we know $\mbox{cd}(G/P_5')=\mbox{cd}(G)\setminus\{p_5\}$, so the subgraph will have vertex set $\{p_1, p_2, p_3, p_4, p_6\}$. If we remove the vertex $p_5$ from the graph, the resulting subgraph is the graph shown to not be possible in \cite{lewis}. The only edges that this subgraph could lose are between two vertices that are adjacent to $p_5$. Thus, we have the following two cases: the edge between $p_4$ and $p_5$, and the edge between $p_3$ and $p_6$. If both edges are lost, we result in a subgraph that violates P\'{a}lfy's condition, and if either one is only lost, we arrive at a subgraph that has five vertices and diameter three, which is not possible by the main theorem of \cite{lewis2} or a graph that violates P\'{a}lfy's condition. Thus, $p_5$ is strongly admissible, and $p_6$ is also. By Lemma \ref{normp} we have that $G$ does not have a normal nonabelian Sylow $p_i$-subgroup for $i\in\{1,4,5,6\}$.
\ecpf\\

Next we show that $G$ does not have a normal Sylow $p_3$-subgroup by showing that the vertex $p_3$ satisfies the hypotheses of Lemma \ref{normiep}.\\

\bcl
The group $G$ does not have a normal Sylow $p_3$-subgroup.
\ecl\\

\bcpf
Consider $p_3$. Then we have that $\pi=\{p_1, p_5,p_6\}$ and $\rho=\{p_2, p_4\}$. We can write $\pi=\pi_1\cup\pi_2$, a disjoint union where $\pi_1=\{p_1\}$ and $\pi_2=\{p_5,p_6\}$. We have shown that $p_1$ is strongly admissible previously, and we know that $p_5$ is adjacent to $p_4$, which is an admissible vertex in $\rho$. Finally, we note that $p_2\in\rho$ and $p_2$ is not adjacent to $p_5$. Thus, the hypotheses of Lemma \ref{normiep} are satisfied, and we have that $G$ does not have a normal Sylow $p_3$-subgroup.
\ecpf\\

We next show that no subgraphs of the graph in Figure \ref{figf} with six vertices can occur.\\

\bcl\label{cl543}
The graph in Figure \ref{figf} has no proper subgraphs with six vertices that occur as a prime character degree graph of a solvable group.
\ecl\\

\bcpf
Suppose $H$ is a group where $\rho(H)=\rho(G)$, $\Delta(H)$ is a proper subgraph of $\Delta(G)$, and $|H|<|G|$. Applying P\'{a}lfy's condition to $p_1$, $p_2$, and $p_6$ in $\Delta(H)$, we see that $p_1$ and $p_2$ must be adjacent in $\Delta(H)$. Again, if we apply P\'alfy's condition to $p_2$, $p_3$, and $p_5$ or $p_6$ in $\Delta(H)$, we see that there must be an edge between $p_3$ and $p_5$ and also between $p_3$ and $p_6$. Applying P\'{a}lfy's condition to $p_2$, $p_5$, and $p_6$, we see that there must be an edge between $p_5$ and $p_6$ in $\Delta(H)$. Again, applying P\'{a}lfy's condition to $p_2$, $p_4$, and $p_3$ in $\Delta(H)$, we see that there must be an edge between $p_2$ and $p_4$. Finally, we claim there must be an edge between $p_1$ and $p_4$ in $\Delta(H)$, since P\'{a}lfy's condition implies that there must be an edge between $p_1$ and $p_4$ or between $p_4$ and $p_6$. If there is an edge between $p_4$ and $p_6$, we arrive at a subgraph with diameter three, which was a graph that was excluded earlier. If we have an edge between $p_1$ and $p_4$, we conclude that this graph along with any possible edge added, can be ruled out by having diameter three, violating P\'{a}lfy's condition, or a graph already ruled out. Thus, $\Delta(H)$ does not occur as the prime character degree graph in any case.
\ecpf\\

\bcl
The group G does not have a normal nonabelian Sylow $p_2$-subgroup.
\ecl\\

\bcpf
We will show that the graph in Figure \ref{figf} satisfies all the conditions of Hypothesis \ref{newhype} for $p=p_2$. Consequently, we have that $\pi =\{p_1, p_4\}$ and $\rho =\{p_3, p_5, p_6\}$.

Since $p_1$ is not adjacent to $p_5$ or $p_6$ and $p_4$ is not adjacent to $p_3$, we satisfy (i) and (ii) of Hypothesis \ref{newhype}. 

We have already shown that $p_1$ and $p_4$ are admissible, and one can check that no connected subgraph with vertex set $\{p\}\cup\{p_1\}\cup\rho$ or $\{p\}\cup\{p_4\}\cup\rho$ occurs as the prime character degree graph for any solvable group $G$. Moreover, the previous claim verifies that no subgraph with vertex set $\{p\}\cup\pi\cup\rho$ occurs either, and so we satisfy (iii).

For the next condition, (iv), we start by investigating the possible subsets $\pi^*\cup\rho^*$, which are: (1) $\{p_1\}\cup\{p_3\}$, (2) $\{p_4\}\cup\{p_5,p_6\}$, (3) $\{p_4\}\cup\{p_5\}$, and (4) $\{p_4\}\cup\{p_6\}$.

Notice that (1) and (2) have their corresponding sets $\tau=\varnothing$, and so there is nothing to verify. For (3), we note that $\tau=\{p_6\}$. By a previous claim, we have verified that $p_6$ is admissible, and one can easily check that no connected subgraph with vertex set $\rho(G)\setminus\{p_6\}=\{p_1,p_2,p_3,p_4,p_5\}$ occurs as $\Delta(G)$ for any solvable group $G$. Observe that a symmetric argument will work for (4), and thus (iv) is satisfied.

The only possible disconnected graph would have components with vertex sets $\{p\}\cup\pi$ and $\rho$, but this is exactly the graph in Figure \ref{fig3and3} which does not occur as the prime character degree graph of any solvable group, and in particular violates P\'alfy's inequality. Thus, (v) is satisfied.

Since we satisfy all the conditions for Hypothesis \ref{newhype}, we apply Theorem \ref{hypethm} and have that $G$ has no normal nonabelian Sylow $p_2$-subgroup.
\ecpf\\

	We have shown that $G$ has no normal nonabelian Sylow subgroups. Let $F$ be the Fitting subgroup of $G$. We note that $\rho(G)=\pi(|G:F|)$, and thus, $\rho(G)=\rho(G/\Phi(G))$ where $\Phi(G)$ is the Frattini subgroup of $G$. Suppose $M$ is a normal subgroup of $G$ so that $\rho(G/M)=\rho(G)$. If $M>1$, then the minimality of $|G|$ would imply that $\Delta(G/M)$ is a proper subgraph of $\Delta(G)$, which violates a previous claim. Thus, we see that if $\rho(G/M)=\rho(G)$, then $M=1$. Now, for $\Phi(G)$, we have that $\rho(G/\Phi(G))=\rho(G)$, so that $\Phi(G)=1$. We may now apply Lemma III 4.4 of \cite{huppert} to see that there is a subgroup $H$ of $G$ so that $G=HF$ and $H\cap F=1$. We let $E$ denote the Fitting subgroup of $H$.\\

\bcl
The Fitting subgroup F of G is a minimal normal subgroup.
\ecl\\

\bcpf
Suppose that there is a normal subgroup $N$ of $G$ so that $1<N<F$. By Theorem III 4.5 of  \cite{huppert}, there is a normal subgroup $M$ of $G$ so that $F=N\times M$. Since $N>1$ and $M>1$, we have $\rho(G/N)\subset\rho(G)$ and $\rho(G/M)\subset\rho(G)$. For any prime $p\in\rho(G)\setminus\rho(G/N)$, we know that $G/N$ has a normal abelian Sylow $p$-subgroup. The class of finite groups with an abelian and normal Sylow $p$-subgroup is a formation, so $p$ must lie in $\rho(G/M)$. Thus, $\rho(G)=\rho(G/N)\cup\rho(G/M)$.\\
\indent If $p\in\rho(G)\setminus\rho(G/N)$, then $p$ is not in $\rho(G/F)=\rho(H)$; so $E$ must then contain the Sylow $p$-subgroup of $H$. Since $p\in\rho(G)$, it follows that $p$ divides $|H|$, and thus $p$ will divide $|E|$. Recall that $\mbox{cd}(G)$ contains a degree divisible by all the prime divisors of $|EF:F|=|E|$. We conclude that $\rho(G)\setminus(\rho(G/M)\cap \rho(G/N))$ lies in a complete subgraph of $\Delta(G)$. Therefore, $\rho(G)\setminus(\rho(G/M)\cap \rho(G/N))$ lies in the subsets: (1) $\{p_1,p_2,p_4\}$, (2) $\{p_1,p_3\}$, (3) $\{p_3,p_5,p_6\}$, or (4) $\{p_4,p_5,p_6\}$.

	Suppose that $(1)$ occurs. This implies that $\{p_3,p_5,p_6\}\subseteq \rho(G/N)\cap\rho(G/M)$. Since $\rho(G/M)\cup\rho(G/N)=\rho(G)$, we know that $\rho(G/M)$ or $\rho(G/N)$ must therefore be one of the following: (i) $\{p_3,p_5,p_6,p_1,p_2\}$, (ii) $\{p_3,p_5,p_6,p_2,p_4\}$, or (iii) $\{p_3,p_5,p_6,p_1,p_4\}$. Assume (i), that is, we have $\rho(G/N)=\{p_1,p_2,p_3,p_5,p_6\}$. By Theorem 5.5 of \cite{lewis2}, $G/N$ has a central Sylow $p_4$-subgroup. However, this would imply $O^{p_4}(G)<G$, a contradiction as we have shown $p_4$ is strongly admissible, which implies $O^{p_4}(G)=G$. A similar argument works for (ii). Finally, assume we have (iii), that is, $\rho(G/N)=\{p_3,p_5,p_6,p_1,p_4\}$. We know $\rho(G/N)\cup\rho(G/M)=\rho(G)$, and so we consider $\rho(G/M)$ in this case. We have that $\{p_3,p_5,p_6,p_2\}\subset\rho(G/M)$. We see that the only possible graph arising from the set of $\rho(G/M)$ must be disconnected and have two components. By Theorem 5.5 of \cite{lewis2}, we have that $G/M$ has either a central Sylow $p_1$-subgroup or a central Sylow $p_4$-subgroup. This implies $O^{p_1}(G)<G$ or $O^{p_4}(G)<G$, respectively. However, this is a contradiction as we showed $p_1$ and $p_4$ are strongly admissible, and hence $O^{p_1}(G)=G=O^{p_4}(G)$.

	Now suppose (2) occurs. This implies that $\{p_2,p_4,p_5,p_6\}\subseteq \rho(G/N)\cap\rho(G/M)$. We consider the possible cases for $\rho(G/N)$ and $\rho(G/M)$. Now we know that $\rho(G/M)$ or $\rho(G/N)$ must contain $p_3$ or $p_1$ as $\rho(G)=\rho(G/N)\cup\rho(G/M)$. Without loss of generality, assume that $\rho(G/N)=\{p_1,p_2,p_4,p_5,p_6\}$ and $\rho(G/M)=\{p_2,p_3,p_4,p_5,p_6\}$. For $\Delta(G/M)$, we have seen that the only possible graph that could occur with this vertex set is the disconnected graph. Again by Theorem 5.5 of \cite{lewis2}, $G/M$ then has a central Sylow $p_1$-subgroup. However, this would imply $O^{p_1}(G)<G$, a contradiction as we have shown $p_1$ is strongly admissible, which implies $O^{p_1}(G)=G$. Thus, (2) cannot occur.

	 Suppose $(3)$ occurs. We then have that $E$ contains a Hall $\{p_3, p_5, p_6\}$-subgroup of $H$. Since cd$(G)$ has a degree divisible by all the primes dividing $|E|$, we have that $|E|$ is divisible by no other primes and $|E|$ is tha Hall $\{p_3, p_5, p_6\}$-subgroup of $H$. Let $\chi\in\text{Irr}(G)$ with $p_3p_5p_6$ dividing $\chi(1)$.  Let $\theta$ be an irreducible constituent of $\chi_{FE}$. Now, $\chi(1)/\theta(1)$ divides $|G:FE|$ and $\chi(1)$ is relatively prime to $|G:FE|$. We determine that $\chi_{FE}=\theta$. Since $p_3,p_5$ and $p_6$ divide $\theta(1)$ and the only possible prime divisors of $a\in\text{cd}(G/FE)$ are $p_1$, $p_2$, or $p_4$, we conclude via Gallagher's theorem that $\text{cd}(G/FE)=\{1\}$ and $G/FE$ is abelian. We now have that $O^{p_4}(G)<G$, a contradiction. So, (3) cannot occur.

	Suppose $(4)$ occurs. Then we have that $\{p_1, p_2,p_3\}\subseteq\rho(G/N)\cap\rho(G/M)$. We know that $\rho(G/M)$ or $\rho(G/N)$ must be one of the following: (i) $\{p_1,p_2,p_3,p_4,p_5\}$, (ii) $\{p_1,p_2,p_3,p_4,p_6\}$, or (iii) $\{p_1,p_2,p_3,p_5,p_6\}$. Note that for each of these cases, the only graph that occurs with the respective vertex set is the disconnected graph. Applying Theorem 5.5 from \cite{lewis2} yields a central Sylow subgroup, which carries us to our contradiction using the admissibility of $p_6$, $p_5$, and $p_4$, respectively. Therefore (4) cannot occur.
\ecpf\\

	Now we apply Lemma \ref{final} for the final contradiction to prove Theorem \ref{th3}. We note that $p_2$ and $p_3$  are adjacent to admissible vertices $p_4$ and $p_5$, $p_3$ is not adjacent to $p_4$, and $p_2$ is not adjacent to $p_5$. Since we have shown $F$ is minimal normal, and the graph we consider has greater than four vertices, Lemma 2.5 of \cite{bissler} is satisfied (see Lemma \ref{final} above), and thus the graph in Figure \ref{figf} is not the prime character degree graph for any solvable group $G$.
\epf

\subsection{Remaining graphs}

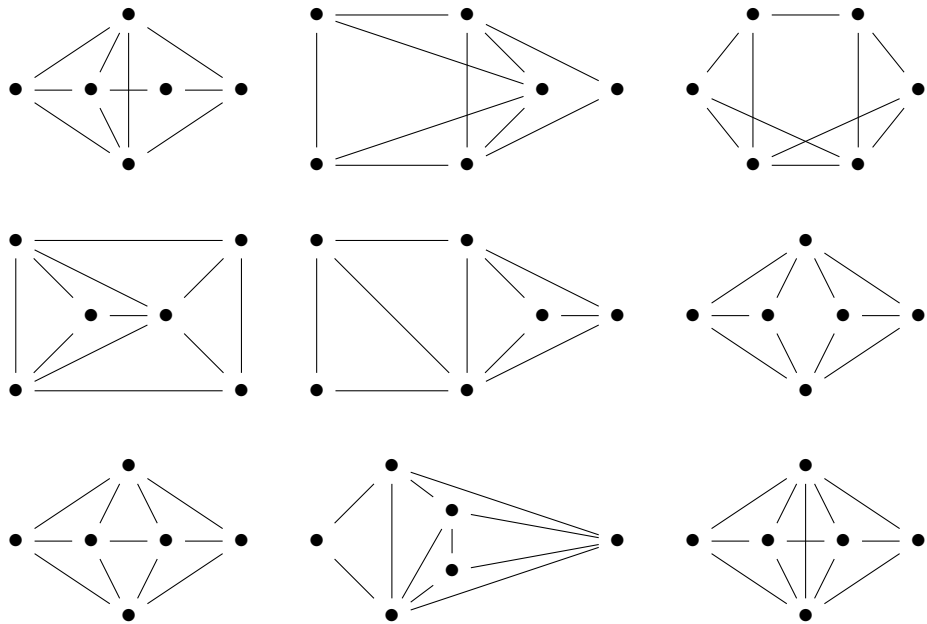
\begin{figure}[h]
    \centering
$
\begin{tikzpicture}[scale=2]
\node (x1) at (0,3.5) {$\bullet$};
\node (x2) at (.5,3.5) {$\bullet$};
\node (x3) at (.75,4) {$\bullet$};
\node (x4) at (.75,3) {$\bullet$};
\node (x5) at (1,3.5) {$\bullet$};
\node (x6) at (1.5,3.5) {$\bullet$};
\path[font=\small,>=angle 90]
(x1) edge node [above] {$ $} (x2)
(x1) edge node [above] {$ $} (x3)
(x1) edge node [above] {$ $} (x4)
(x2) edge node [above] {$ $} (x3)
(x2) edge node [above] {$ $} (x4)
(x2) edge node [above] {$ $} (x5)
(x3) edge node [above] {$ $} (x4)
(x3) edge node [above] {$ $} (x6)
(x4) edge node [above] {$ $} (x6)
(x5) edge node [above] {$ $} (x6);
\node (y1) at (2,4) {$\bullet$};
\node (y2) at (2,3) {$\bullet$};
\node (y3) at (3,4) {$\bullet$};
\node (y4) at (3,3) {$\bullet$};
\node (y5) at (3.5,3.5) {$\bullet$};
\node (y6) at (4,3.5) {$\bullet$};
\path[font=\small,>=angle 90]
(y1) edge node [above] {$ $} (y2)
(y1) edge node [above] {$ $} (y3)
(y1) edge node [above] {$ $} (y5)
(y2) edge node [above] {$ $} (y4)
(y2) edge node [above] {$ $} (y5)
(y3) edge node [above] {$ $} (y4)
(y3) edge node [above] {$ $} (y5)
(y3) edge node [above] {$ $} (y6)
(y4) edge node [above] {$ $} (y5)
(y4) edge node [above] {$ $} (y6);
\node (z1) at (4.5,3.5) {$\bullet$};
\node (z2) at (4.9,4) {$\bullet$};
\node (z3) at (4.9,3) {$\bullet$};
\node (z4) at (5.6,4) {$\bullet$};
\node (z5) at (5.6,3) {$\bullet$};
\node (z6) at (6,3.5) {$\bullet$};
\path[font=\small,>=angle 90]
(z1) edge node [above] {$ $} (z2)
(z1) edge node [above] {$ $} (z3)
(z1) edge node [above] {$ $} (z5)
(z2) edge node [above] {$ $} (z3)
(z2) edge node [above] {$ $} (z4)
(z3) edge node [above] {$ $} (z5)
(z3) edge node [above] {$ $} (z6)
(z4) edge node [above] {$ $} (z5)
(z4) edge node [above] {$ $} (z6)
(z5) edge node [above] {$ $} (z6);
\node (a1) at (0,2.5) {$\bullet$};
\node (a2) at (0,1.5) {$\bullet$};
\node (a3) at (.5,2) {$\bullet$};
\node (a4) at (1,2) {$\bullet$};
\node (a5) at (1.5,2.5) {$\bullet$};
\node (a6) at (1.5,1.5) {$\bullet$};
\path[font=\small,>=angle 90]
(a1) edge node [above] {$ $} (a2)
(a1) edge node [above] {$ $} (a3)
(a1) edge node [above] {$ $} (a4)
(a1) edge node [above] {$ $} (a5)
(a2) edge node [above] {$ $} (a3)
(a2) edge node [right] {$ $} (a4)
(a2) edge node [above] {$ $} (a6)
(a3) edge node [above] {$ $} (a4)
(a4) edge node [above] {$ $} (a5)
(a4) edge node [above] {$ $} (a6)
(a5) edge node [right] {$ $} (a6);
\node (b1) at (2,2.5) {$\bullet$};
\node (b2) at (2,1.5) {$\bullet$};
\node (b3) at (3,2.5) {$\bullet$};
\node (b4) at (3,1.5) {$\bullet$};
\node (b5) at (3.5,2) {$\bullet$};
\node (b6) at (4,2) {$\bullet$};
\path[font=\small,>=angle 90]
(b1) edge node [above] {$ $} (b2)
(b1) edge node [above] {$ $} (b3)
(b1) edge node [above] {$ $} (b4)
(b2) edge node [above] {$ $} (b4)
(b3) edge node [right] {$ $} (b4)
(b3) edge node [above] {$ $} (b5)
(b3) edge node [above] {$ $} (b6)
(b4) edge node [above] {$ $} (b5)
(b4) edge node [above] {$ $} (b6)
(b5) edge node [right] {$ $} (b6);
\node (c1) at (4.5,2) {$\bullet$};
\node (c2) at (5.25,2.5) {$\bullet$};
\node (c3) at (5,2) {$\bullet$};
\node (c4) at (5.25,1.5) {$\bullet$};
\node (c5) at (5.5,2) {$\bullet$};
\node (c6) at (6,2) {$\bullet$};
\path[font=\small,>=angle 90]
(c1) edge node [right] {$ $} (c2)
(c1) edge node [above] {$ $} (c3)
(c1) edge node [above] {$ $} (c4)
(c2) edge node [above] {$ $} (c3)
(c2) edge node [right] {$ $} (c5)
(c2) edge node [above] {$ $} (c6)
(c3) edge node [above] {$ $} (c4)
(c4) edge node [above] {$ $} (c5)
(c4) edge node [above] {$ $} (c6)
(c5) edge node [right] {$ $} (c6);
\node (d1) at (0,.5) {$\bullet$};
\node (d2) at (.75,1) {$\bullet$};
\node (d3) at (.5,.5) {$\bullet$};
\node (d4) at (.75,0) {$\bullet$};
\node (d5) at (1,.5) {$\bullet$};
\node (d6) at (1.5,.5) {$\bullet$};
\path[font=\small,>=angle 90]
(d1) edge node [right] {$ $} (d2)
(d1) edge node [above] {$ $} (d3)
(d1) edge node [above] {$ $} (d4)
(d2) edge node [above] {$ $} (d3)
(d2) edge node [right] {$ $} (d5)
(d2) edge node [above] {$ $} (d6)
(d3) edge node [above] {$ $} (d4)
(d3) edge node [above] {$ $} (d5)
(d4) edge node [above] {$ $} (d5)
(d4) edge node [above] {$ $} (d6)
(d5) edge node [right] {$ $} (d6);
\node (e1) at (2,.5) {$\bullet$};
\node (e2) at (2.5,1) {$\bullet$};
\node (e3) at (2.9,.7) {$\bullet$};
\node (e4) at (2.5,0) {$\bullet$};
\node (e5) at (2.9,.3) {$\bullet$};
\node (e6) at (4,.5) {$\bullet$};
\path[font=\small,>=angle 90]
(e1) edge node [right] {$ $} (e2)
(e1) edge node [above] {$ $} (e4)
(e2) edge node [above] {$ $} (e3)
(e2) edge node [above] {$ $} (e4)
(e2) edge node [above] {$ $} (e6)
(e3) edge node [above] {$ $} (e4)
(e3) edge node [above] {$ $} (e5)
(e3) edge node [above] {$ $} (e6)
(e4) edge node [above] {$ $} (e5)
(e4) edge node [above] {$ $} (e6)
(e5) edge node [right] {$ $} (e6);
\node (f1) at (4.5,.5) {$\bullet$};
\node (f2) at (5.25,1) {$\bullet$};
\node (f3) at (5,.5) {$\bullet$};
\node (f4) at (5.25,0) {$\bullet$};
\node (f5) at (5.5,.5) {$\bullet$};
\node (f6) at (6,.5) {$\bullet$};
\path[font=\small,>=angle 90]
(f1) edge node [right] {$ $} (f2)
(f1) edge node [above] {$ $} (f3)
(f1) edge node [above] {$ $} (f4)
(f2) edge node [above] {$ $} (f3)
(f2) edge node [above] {$ $} (f4)
(f2) edge node [right] {$ $} (f5)
(f2) edge node [above] {$ $} (f6)
(f3) edge node [above] {$ $} (f4)
(f3) edge node [above] {$ $} (f5)
(f4) edge node [above] {$ $} (f5)
(f4) edge node [above] {$ $} (f6)
(f5) edge node [right] {$ $} (f6);
\end{tikzpicture}
$
    \caption{Unknown six-vertex graphs}
    \label{figunknown}
\end{figure}

There are still nine remaining graphs with six vertices that have yet to be classified (see Figure \ref{figunknown}).

\end{document}